\documentclass{SciPost}

\hypersetup{
    colorlinks,
    linkcolor={red!50!black},
    citecolor={blue!50!black},
    urlcolor={blue!80!black}
}

\usepackage[bitstream-charter]{mathdesign}
\DeclareSymbolFont{usualmathcal}{OMS}{cmsy}{m}{n}
\DeclareSymbolFontAlphabet{\mathcal}{usualmathcal}

\fancypagestyle{SPstyle}{
\fancyhf{}
\lhead{\colorbox{scipostblue}{\bf \color{white} ~SciPost Physics Proceedings }}
\rhead{{\bf \color{scipostdeepblue} ~Submission }}

\fancyfoot[C]{\textbf{\thepage}}
}

\usepackage{amsthm,amscd,euscript}
\usepackage{algorithm,algpseudocode}
\usepackage{tikz}
\usepackage{booktabs}
\usepackage{float}
\usepackage{array}
\usepackage{multirow}
\usepackage{seqsplit}
\providecommand{\seqnum}[1]{#1}

\definecolor{up}{HTML}{378ADD}\definecolor{dn}{HTML}{D85A30}
\definecolor{upbg}{HTML}{B5D4F4}\definecolor{dnbg}{HTML}{F5C4B3}
\definecolor{upd}{HTML}{0C447C}\definecolor{dnd}{HTML}{712B13}

\theoremstyle{plain}

\theoremstyle{definition}
\newtheorem{theor}{Theorem}[section]
\newtheorem{lemma}[theor]{Lemma}

\newtheorem{conjecture}[theor]{Conjecture}

\theoremstyle{definition}
\newtheorem{claim}[theor]{Claim}

\newtheorem*{convention}{Convention}

\newtheorem{prop}[theor]{Proposition}

\newtheorem{define}[theor]{Definition}

\newtheorem{notation}[theor]{Notation}
\newtheorem{ex}[theor]{Example}
\theoremstyle{remark}
\newtheorem{rem}[theor]{Remark}

\newcommand{\BBR}{\mathbb{R}}
\newcommand{\BBN}{\mathbb{N}} 
\newcommand{\BBS}{\mathbb{S}}
\newcommand{\BBZ}{\mathbb{Z}}

\newcommand{\dd}{\partial}

\DeclareMathOperator{\id}{id}
\DeclareMathOperator{\Alt}{Alt}
\DeclareMathOperator{\Wronsk}{Wronsk}

\DeclareMathOperator{\consta}{\boldsymbol{c}}

\newcommand{\by}[1]{\textit{{#1}}}

\providecommand{\nopindent}{\noindent}
\begin{document}

\pagestyle{SPstyle}

\begin{center}{\Large \textbf{\color{scipostdeepblue}{
New integer sequence OEIS\:A392\,714 counts Wronskians, II.\\
Computing within the permutation group algebra
}}}\end{center}

\begin{center}\textbf{
Kian C.\ Shah\textsuperscript{1} and
Arthemy V.\ Kiselev\textsuperscript{1$\star$}
}\end{center}

\begin{center}
{\bf 1} Bernoulli Institute for Mathematics, Computer Science and Artificial Intelligence, University of Groningen, Groningen, The Netherlands
\\[\baselineskip]
$\star$ \href{mailto:a.v.kiselev@rug.nl}{\small a.v.kiselev@rug.nl}
\end{center}

\definecolor{palegray}{gray}{0.95}
\begin{center}
\colorbox{palegray}{
  \begin{tabular}{ll}
 \hspace{-11pt} \begin{minipage}{0.18\textwidth}
    \includegraphics[width=28mm]{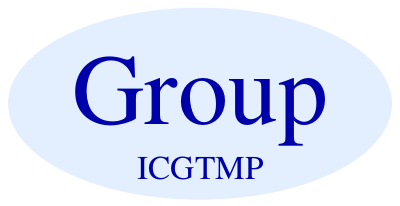}
  \end{minipage}
  &
  \begin{minipage}{0.75\textwidth}
    \begin{center}
    {\it 36th International Colloquium on Group Theoretical Methods in Physics}\\
    {\it Valladolid, 13-17 July 2026} \\
    \doi{10.21468/SciPostPhysProc.?}\\
    \end{center}
  \end{minipage}
\end{tabular}
}
\end{center}

\section*{\color{scipostdeepblue}{Abstract}}
\textbf{\boldmath{In the context of strongly homotopy Lie algebras, we study Wronskians as $N$-ary Lie brackets. Namely, let $w_j(x)\cdot\partial_x^p$ be differential operators with smooth coefficients and of strict order $p\in\BBN_{\geqslant 1}$ on the real line~$\BBR$; by alternating $N=2p$ of them, we get $\Alt(w_1\partial_x^p,\ldots,w_N\partial_x^p)=\consta(p)\cdot\Wronsk(\{w_j\})\cdot\partial_x^p$. Every $\sigma\in\BBS_{N(=2p)}$ contributes with $(-)^\sigma\cdot(w_{\sigma(1)}\partial_x^p\circ\cdots\circ w_{\sigma(N-1)}\partial_x^p)(w_{\sigma(N)})$ to the r.h.s.; whenever we put $w_j:=x^{j-1}$, this is $(-)^\sigma c(\sigma,p)\in\BBZ$, and the support of the map $\sigma\mapsto(-)^\sigma c(\sigma,p)$ is the set~$\Phi_p$ of late\/-\/growing permutations. For instance, $\Phi_{p=1}=\{\id\}$ and Lie brackets fix $\consta(p=1)=1$. We discover that computing $\consta(p>1)$ can be done fast in the permutation group algebra (indeed, by subtracting $\{\id-\sigma\}$ termwise). This algebraic approach not only gives us the exact values up to $p=7$ but also provides the lower bound $v_p(\consta(p))\geqslant p-1$ for the $p$-adic valuations of~$\consta(p)$ at any prime~$p$.}}

\vspace{\baselineskip}

\noindent\textcolor{white!90!black}{%
\fbox{\parbox{0.975\linewidth}{%
\textcolor{white!40!black}{\begin{tabular}{lr}%
  \begin{minipage}{0.6\textwidth}%
    {\small Copyright attribution to authors. \newline
    This work is a submission to SciPost Phys. Proc. \newline
    License information to appear upon publication. \newline
    Publication information to appear upon publication.}
  \end{minipage} & \begin{minipage}{0.4\textwidth}
    {\small Received Date \newline Accepted Date \newline Published Date}%
  \end{minipage}
\end{tabular}}
}}
}

\vspace{10pt}
\noindent\rule{\textwidth}{1pt}
\tableofcontents
\noindent\rule{\textwidth}{1pt}
\vspace{10pt}

\section{Introduction}\label{sec:intro}

There is a natural way to construct $N$-ary Lie brackets that satisfy the table of Jacobi identities for strongly homotopy Lie algebras (see \cite{Dzhuma2014,LadaStasheff} and \cite{ForKac,PRG25}): namely, by alternating $N=2p$ differential operators of strict order $p\in\BBN_{\geqslant 1}$ on the line $\BBR$. The resulting bracket is equal to the product of the universal constant $\consta(p)\in\BBZ$ and Wronskian determinant of the $N$ old coefficients; the new operator is again of strict order~$p$. This is a natural generalization of the commutator $[\cdot,\cdot]$ for smooth vector fields on $\BBR$ at $p=1$.

In this paper we describe an algorithm to compute the constants $\consta(p\geqslant 4)$ beyond the base $\consta(1)=1$ from Lie theory and easily obtained $\consta(2)=2$ and $\consta(3)=90$; we discover that the integer sequence $\{\consta(p)\}_{p\in\BBN}$ is new \cite{oeisA392714}. The algorithm, which yields all the exact values up to $\consta(p=7)\approx 8.587\cdot 10^{34}$, deploys calculation within the permutation group algebra; the same technique, referring to the deviations of permutations at hand from the identity mapping, helps us to prove the lower bound $v_p(\consta(p))\geqslant p-1$ for the power of primes~$p$ in the $p$-adic valuations of the constants $\consta(p)$. The algorithm which we presently formulate can be accelerated at the price of memoization (see \cite{Method}, where $\consta(p=18)\approx 4.881\cdot 10^{462}$ is attained; the values known from this text (before the speedup) were used for cross-verification of the larger output).

We note also that the late-growing permutations (see footnote~\ref{FootDefLateGrow} on p.~\pageref{FootDefLateGrow} for definition) are the set which survives from $\BBS_{2p}$ through the filter whenever the $p$th order differential operators' coefficients are chosen to be the monomials $1,\ldots,x^{N-1}$, so that it is their degrees $\geqslant 0$ which select the permutations of $\BBS_{2p}$ that actually contribute to $\consta(p)$. In \cite{Method,cpWron} we detect that the sets of late-growing permutations are almost equally partitioned by the parity-even vs odd components; this we prove in \cite{OddEven}. The fact of balance itself reveals a hidden asymmetry of the contributing permutations: despite the almost-equal number of positive and negative contributions to $\consta(p)\in\BBZ$ in the course of alternating summation, the values $\consta(p)$ are always positive and grow as fast as $\consta(p)\sim p^{2p^2}$ in the limit $p\to+\infty$. The nature of this asymmetry is yet to be explored.

\begin{center}
\leavevmode\vbox{\hrule width 0.45\linewidth \vskip 2pt \hrule width 0.45\linewidth}
\end{center}

\noindent This paper is structured as follows. We first recall in \S\ref{sec:prelim} that by alternating $N=2p$ differential operators of strict order $p\in\BBN_{\geqslant 1}$ on the real line $\BBR$, one obtains the operator of the same strict order~$p$; the new coefficient is the Wronskian determinant of the $N$ old coefficients times the universal constant $\consta(p)$. The problem is to find its values for as many arguments $p>1$ as possible. By postulating that the weights be monomials $1,\ldots,x^{N-1}$, we convert the problem of analysis into one from combinatorics in \S\ref{sec:monomials}: it is the monomial degrees ($\geqslant 0$) which we keep track of; every permutation from the alternating sum contributes by the product of falling factorials to the constant near Wronskian (if it does at all, whence the filter of late-growing permutations in \S\ref{sec:algorithm}). As soon as the formula of $\consta(p)$ is obtained (see Eq.~\eqref{eq:closedform} in \S\ref{sec:closedform}, illustrated in Example~\ref{ex:peq2} on p.~\pageref{ex:peq2}), we explain in \S\ref{sec:deviation} the key idea of the originally developed fast algorithm to reach the values $\consta(p\leqslant 7)$: here we work in the permutation group algebra, literally taking the termwise differences of permutations. The values $\consta(p\leqslant 7)$ are listed in \S\ref{sec:results}; we discover that the positive integers $\{\consta(p)\}_{p\in\BBN}$, originally introduced in~\cite{Dzhuma2014}, form a new sequence \cite{oeisA392714}; its growth rate is $\consta(p)\sim p^{2p^2}$ as $p\to+\infty$ (O.\,Zaboronski, private communication).\footnote{\label{FootZabor}
\by{O.\,Zaboronski} (June 2026) $\consta(p)$ as the solution to a counting problem.
IH\'ES, 7~p. (work in progress), unpublished.} We study the prime factorisations of the known $\consta(p)$; in \S\ref{sec:lowerbound} we prove that for prime indices~$p$, the factorization contains $p$ to the power $\geqslant p-1$; given the evidence, we conjecture that the lower bound is actually an equality.

\section{Preliminaries}\label{sec:prelim}

\begin{define}\label{def:altcomp}
Let $p \in \mathbb{Z}_{\geqslant 1}$, and set $N = 2p$; let $w_1$,\ $\ldots$,\ $w_N$
be the coefficients of differential operators of strict order $p$ on an
interval $I \subseteq \mathbb{R}$ with an affine coordinate~$x$. For a smooth
function $f \in C^\infty(\mathbb{R})$, the \emph{alternating
composition} of these differential operators is
\begin{equation}\label{eq:opdef}
\mathcal{A}_p[w_1,\ldots,w_N](f)(x) \mathrel{{:}{=}}
\sum_{\sigma \in \BBS_N} (-1)^\sigma\,
w_{\sigma(1)}\partial_x^p \circ w_{\sigma(2)}\partial_x^p
\circ \cdots \circ w_{\sigma(N)}\partial_x^p\,(f(x));
\end{equation}
the sum runs over all permutations $\sigma$ of the set
$\{1,\ldots,N\}$ and $(-1)^\sigma$ denotes the sign of~$\sigma$.
\end{define}

\begin{convention}
Throughout this paper, we adopt the one-line notation for permutations, in
which $\sigma = (\sigma(1), \sigma(2), \ldots, \sigma(N))$ records the
image at each position. That is, $\sigma(k)$ denotes the index of the
object occupying position $k$ after the permutation is applied. Under this
convention, the identity element of $\BBS_N$ is represented as
$\sigma = (1, 2, \ldots, N)$. For instance, in $\BBS_4$, the identity
$\sigma = (1,2,3,4)$ corresponds to the map
$(1,\, x,\, x^2,\, x^3) \;\longmapsto\; (1,\, x,\, x^2,\, x^3)$,
while $\sigma = (4,1,2,3)$ corresponds to
$(1,\, x,\, x^2,\, x^3) \;\longmapsto\; (x^3,\, 1,\, x,\, x^2)$.
\end{convention}

\begin{theor}[{\cite{ForKac,PRG25} and~\cite{Dzhuma2014}}]\label{thm:main}
For any smooth function $f$ and all weight functions $w_1, \ldots, w_N$,
\begin{equation}\label{EqMain}
\mathcal{A}_p[w_1,\ldots,w_N](f)(x)
= \consta(p) \cdot
\Wronsk \bigl(w_1, \ldots, w_N\bigr)(x) \cdot \partial_x^p(f(x)),
\end{equation}
where $\consta(p)$ is a universal constant, independent of $f$ and
of the weight functions.
\end{theor}

The proof of this theorem is outlined in~\cite{Dzhuma2014} and~\cite{cpWron,Method}.

\begin{lemma}\label{lem:order}
Put $f(x)=x^p/p!$ in~\eqref{EqMain}, so that $\partial_x^p(f)=1$. The rightmost operator $w_{\sigma(N)}\partial_x^p$ in every term of~\eqref{eq:opdef} takes this~$f$ to the weight~$w_{\sigma(N)}$, and \eqref{EqMain} reduces to
the following equation:
\begin{equation}\label{eq:reduced}
\sum_{\sigma \in \BBS_N} (-1)^\sigma\, \bigl(w_{\sigma(1)}\partial_x^p \circ
\cdots \circ w_{\sigma(N-1)}\partial_x^p\bigr)\bigl(w_{\sigma(N)}\bigr)
\;=\; \consta(p) \cdot \Wronsk(w_1, \ldots, w_N).
\end{equation}
\end{lemma}

\begin{define}\label{def:groupalg}
The \emph{group algebra} $\BBZ[\BBS_N]$ of the permutation group~$\BBS_N$ is the ring of formal sums $\sum_{\sigma\in\BBS_N} n_\sigma\cdot\sigma$ with integer coefficients~$n_\sigma$; these sums are added coefficientwise and multiplied by extending the composition of permutations bilinearly. The \emph{alternator} is the element
\[
\mathbf{a}_N := \sum_{\sigma\in\BBS_N} (-1)^\sigma\cdot\sigma \ \in\ \BBZ[\BBS_N].
\]
\end{define}

\begin{define}\label{def:ev}
Let $w=(w_1,\ldots,w_N)$ be a tuple of weights. The \emph{evaluation map} $\mathrm{ev}_w$ is the $\BBZ$-linear map from $\BBZ[\BBS_N]$ to the associative algebra of differential operators with smooth coefficients such that, for every $\sigma\in\BBS_N$,
\[
\mathrm{ev}_w(\sigma) := w_{\sigma(1)}\partial_x^p\circ w_{\sigma(2)}\partial_x^p\circ\cdots\circ w_{\sigma(N)}\partial_x^p .
\]
\end{define}

\begin{rem}\label{rem:evalt}
By construction, the right-hand side of~\eqref{eq:opdef} is $\mathrm{ev}_w(\mathbf{a}_N)(f)$. Theorem~\ref{thm:main} is the equality of differential operators,
\begin{equation}\label{eq:evalt}
\mathrm{ev}_w(\mathbf{a}_N) = \consta(p)\cdot\Wronsk(w_1,\ldots,w_N)\cdot\partial_x^p .
\end{equation}
\end{rem}

The problem which we study in this paper is to find the values of the constant $\consta(p)$ in equality~\eqref{eq:evalt} by calculating within the group algebra $\BBZ[\BBS_N]$.

For small values $p \leqslant 3$ we verify identity \eqref{EqMain} directly: by hand for $p = 1$, whence $\consta(1) = 1$, and for $p = 2$, whence we get
$\consta(2) = 2$ (cf.\ Example~\ref{ex:peq2} below). For $p = 3$,
using the package \textsf{Jets} for differential calculus
(see~\cite{BaranMarvan_Jets,Marvan2009}) in \textsf{Maple}, we
expand the left\/-\/hand side of \eqref{EqMain} by taking the entire sum over
$6! = 720$ permutations (here $N = 2 \cdot 3 = 6$ and all differential
operators have order $p = 3$); we find $\consta(3) = 90$. The values
of $\consta(p)$ for small $p$ are listed in~\eqref{eq:smallp}:
\begin{equation}\label{eq:smallp}
\underset{N!=2}{\consta(1)=1}\ \text{(Lie)};\qquad
\underset{N!=24}{\consta(2)=2};\qquad
\underset{N!=720}{\consta(3)=90}.
\end{equation}

Dzhumadil'daev provided a formula of $\consta(p)$ and calculated three more values, namely $\consta(p\leqslant 6)$, in~\cite[\S5]{Dzhuma2014}. To learn more values $\consta(p\geqslant 7)$, we now use a strategy different from the one in~\cite{Dzhuma2014}, directly referring to the definition of $\consta(p)$ by Eq.~\eqref{EqMain}.

\section{Combinatorics of the monomial case}\label{sec:monomials}

\noindent%
To make the right\/-\/hand side of~\eqref{eq:reduced} shorter and extract the
constants $\consta(p)$ explicitly, we choose specific weight
functions, namely, let them be monomials:
\begin{equation}\label{EqMonomWeight}
w_n(x) = x^{n-1}, \qquad n = 1, \ldots, N.
\end{equation}

\begin{rem}
Note that $\deg(w_{\sigma(k)}) = \sigma(k) - 1$ whenever $w_k = x^{k-1}$.
\end{rem}

\begin{lemma}\label{lem:wronskian}
The Wronskian of $1, x, x^2, \ldots, x^{N-1}$ equals
\[
\Wronsk(1,\, x,\, x^2,\, \ldots,\, x^{N-1})
\;=\; \prod_{k=0}^{N-1} k!\,.
\]
\end{lemma}

\begin{proof}
The Wronskian matrix of $1, x, \ldots, x^{N-1}$ is upper
triangular, with diagonal entries $0!$,\ $1!$,\ $2!$,\ $\ldots$,\ $(N-1)!$. Thus, the
determinant is the product of these entries.
\end{proof}

To understand the action of the iterated operators in \eqref{eq:opdef} in
the monomial case, consider two monomials $x^\alpha$ and $x^\beta$ with
$\alpha, \beta \in \mathbb{Z}_{\geqslant 0}$. The operator $x^\beta\partial_x^p$
acts on $x^\alpha$: whenever $\alpha < p$ the result is zero, otherwise
\begin{equation}\label{eq:opaction}
x^\beta\,\partial_x^p\,(x^\alpha)
\;=\; x^\beta\cdot\frac{\alpha!}{(\alpha-p)!}\, x^{\alpha-p}
\;=\; \frac{\alpha!}{(\alpha-p)!}\, x^{\alpha+\beta-p}.
\end{equation}
\begin{notation}\label{not:falling}
For integers $\alpha\geqslant p$, the \textsl{falling factorial} is
\[
\alpha^{\underline{p}}
= \alpha!/(\alpha-p)!
= \alpha(\alpha-1)\cdots(\alpha-p+1).
\]
\end{notation}

\noindent Each application of the weighted operator $x^\beta\partial_x^p$, when
nonzero in the output, multiplies the argument $x^\alpha$ by the falling
factorial~$\alpha^{\underline{p}}$; also, the operator $x^\beta\partial_x^p$ shifts the exponent by
$\alpha \mapsto \alpha + \beta - p$.

\begin{rem}
The vanishing condition $\alpha < p$ immediately implies that not every
permutation $\sigma \in \BBS_N$ contributes to the alternating sum in~\eqref{eq:reduced}. Any permutation for which the innermost weight function
$w_{\sigma(N)} = x^{\sigma(N)-1}$ satisfies $\sigma(N) - 1 < p$ yields zero
contribution; clearly, such a permutation may be discarded at the outset.
\end{rem}

\begin{ex}\label{ex:filter}
Let $p = 2$ and $N = 4$, so that the weight functions are
$\{1, x, x^2, x^3\}$. Consider the permutation $\sigma = (3,1,4,2)$,
corresponding to
$(1,\, x,\, x^2,\, x^3) \mapsto (x^2,\, 1,\, x^3,\, x)$,
so that $\sigma(4) = 2$ and the innermost operator acts as
$\partial_x^2(w_2) = \partial_x^2(x) = 0$. The entire term arising from
this permutation thus vanishes, and this $\sigma$ does not contribute to the
sum. More generally, any permutation with $\sigma(N) \leqslant p$ is eliminated
by this filter; for $p = 2$, $N = 4$ this removes all $\sigma$ with
$\sigma(4) \in \{1, 2\}$, reducing the number of permutations to be
evaluated from $4! = 24$ to those with $\sigma(4) \in \{3, 4\}$.
\end{ex}

As illustrated by the above example, we can now shift our focus
from the monomial weight functions to the permutations themselves; from now on, we
work with permutations directly.

\begin{define}\label{def:csigma}
Let the weights be the monomials~\eqref{EqMonomWeight}. For $\sigma\in\BBS_N$, denote by
\[
c(\sigma,p) := \bigl(w_{\sigma(1)}\partial_x^p\circ w_{\sigma(2)}\partial_x^p\circ\cdots\circ w_{\sigma(N-1)}\partial_x^p\bigr)(w_{\sigma(N)})
\]
the term which $\sigma$ contributes, without its sign $(-1)^\sigma$, to the left-hand side of~\eqref{eq:reduced}.
\end{define}

\begin{lemma}\label{lem:calt}
\textup{(}\textit{i}\textup{)} For every $\sigma\in\BBS_N$, the function $c(\sigma,p)$ is a constant, and this constant is a nonnegative integer.

\textup{(}\textit{ii}\textup{)} Let $c(\,\cdot\,,p)\colon\BBZ[\BBS_N]\to\BBZ$ be the $\BBZ$-linear map which takes every $\sigma\in\BBS_N$ to $c(\sigma,p)$. Then
\begin{equation}\label{eq:calt}
c(\mathbf{a}_N,p) = \consta(p)\cdot\prod_{k=0}^{N-1}k!\,.
\end{equation}
\end{lemma}

\begin{proof}
(\textit{i}) By~\eqref{eq:opaction}, $c(\sigma,p)$ is either zero or the product of $N-1$ falling factorials, which are positive integers, and of a monomial in~$x$ of degree
\[
\sum_{i=1}^{N}\bigl(\sigma(i)-1\bigr)-(N-1)\,p = \tfrac{1}{2}N(N-1)-(N-1)\,p = 0,
\]
because $N=2p$.

(\textit{ii}) For the monomial weights, $c(\mathbf{a}_N,p)$ is the left-hand side of~\eqref{eq:reduced}; its right-hand side is $\consta(p)\cdot\prod_{k=0}^{N-1}k!$ by Lemma~\ref{lem:wronskian}.
\end{proof}

Thus the constant $\consta(p)$ is found by evaluating one $\BBZ$-linear map on one element, the alternator, of the group algebra $\BBZ[\BBS_N]$.

\begin{lemma}\label{lem:sigma1}
For every contributing permutation $\sigma \in \BBS_N$, we have
$\sigma(1) = 1$\textup{;} that is, the weight function $w_{\sigma(1)} = w_1 = 1$ is
necessarily fixed in the leading leftmost position.
\end{lemma}

\begin{proof}
There are two kinds of permutations: the ones which keep the leftmost
monomial $w_1 = 1$ in its first place, and those which send it away, so
that a positive-degree monomial $x^{\alpha}$ becomes the leftmost. Suppose
the latter, so that it is this positive-degree monomial $x^{\alpha}$ which
itself is not differentiated, but is multiplied by a nonnegative-degree
monomial resulting from all the differentiations which stood to the right
of it. This yields a contribution of degree $\geqslant \alpha$ to the sum. But
we know from Lemma~\ref{lem:calt}(\textit{i}) that every term $c(\sigma,p)$ in the left-hand side of \eqref{eq:reduced} is a real number, i.e.\ it is of polynomial degree
zero. Therefore, every permutation sending $1$ away cannot contribute to
the left-hand side of \eqref{eq:reduced}.
\end{proof}

\section{Algorithm to filter permutations and a faster generative procedure}
\label{sec:algorithm}

\noindent\textbf{Generator (ver. 1).} We start by generating the set of all $N!$ permutations of ${1, 2, \dots, N}$ by a standard recursive construction. Then, any permutation with $\sigma(1) \neq 1$ is discarded.
The remaining permutations are tested by simulating the exponent-tracking process from right to left.
\begin{define}\label{def:Tk}
For $\sigma\in\BBS_N$, the partial sums are
\[
T_k \;:=\; \sum_{i=N-k+1}^{N} \bigl(\sigma(i) - 1-p\bigr), \qquad k = 1, \dots, N-1.
\]
\end{define}

\noindent The quantity $T_k$ tracks whether the exponent of $x$ remains non-negative after each successive operator application. If $T_k < 0$ for some $k$, the exponent would become negative at that stage, causing the term to vanish; this permutation is therefore rejected. Only those permutations for which $T_k \geqslant 0$ for all $k$ are retained.
\begin{notation}
We denote by \(\Phi_p\) the set of all the remaining permutations (from \(\BBS_N\)) which pass the filter and contribute to the sum in the l.-h.s.\ of \eqref{eq:reduced}; note \(\Phi_p \subsetneq \BBS_N\).
\end{notation}

\begin{rem}\label{rem:support}
In terms of Definition~\ref{def:csigma}, the set $\Phi_p$ consists of those $\sigma\in\BBS_N$ for which $c(\sigma,p)\neq 0$. Thus
\[
c(\mathbf{a}_N,p) = \sum_{\sigma\in\Phi_p}(-1)^\sigma\,c(\sigma,p).
\]
\end{rem}

\noindent\textbf{Generator (ver. 2).}The same idea as above can be used to make a generative procedure, which avoids generating invalid permutations in the first place. Rather than constructing all $N!$ permutations \textit{ab initio} and filtering \textit{post hoc}, it builds permutations incrementally via backtracking, enforcing the non-negativity condition $T_k \geqslant 0$ at each step of the construction.

\begin{figure}[htb]
\centering
\begin{tikzpicture}[level distance=11mm, level 1/.style={sibling distance=42mm}, level 2/.style={sibling distance=21mm}, every node/.style={font=\small}]
\node {$\sigma(1)=1$}
 child {node {$\sigma(4)=2$\ \ $\times$}}
 child {node {$\sigma(4)=3$}
   child {node {$\sigma(3)=2$\ \ $\times$}}
   child {node {$\sigma(3)=4$} child {node {$(1,2,4,3)$}}}}
 child {node {$\sigma(4)=4$}
   child {node {$\sigma(3)=2$} child {node {$(1,3,2,4)$}}}
   child {node {$\sigma(3)=3$} child {node {$(1,2,3,4)$}}}};
\end{tikzpicture}
\caption{The generator (ver.~2) for $p=2$, $N=4$: the positions are filled from right to left, and a branch is pruned ($\times$) as soon as $T_k<0$. The three leaves are the permutations in~$\Phi_2$.}\label{fig:pruning}
\end{figure}
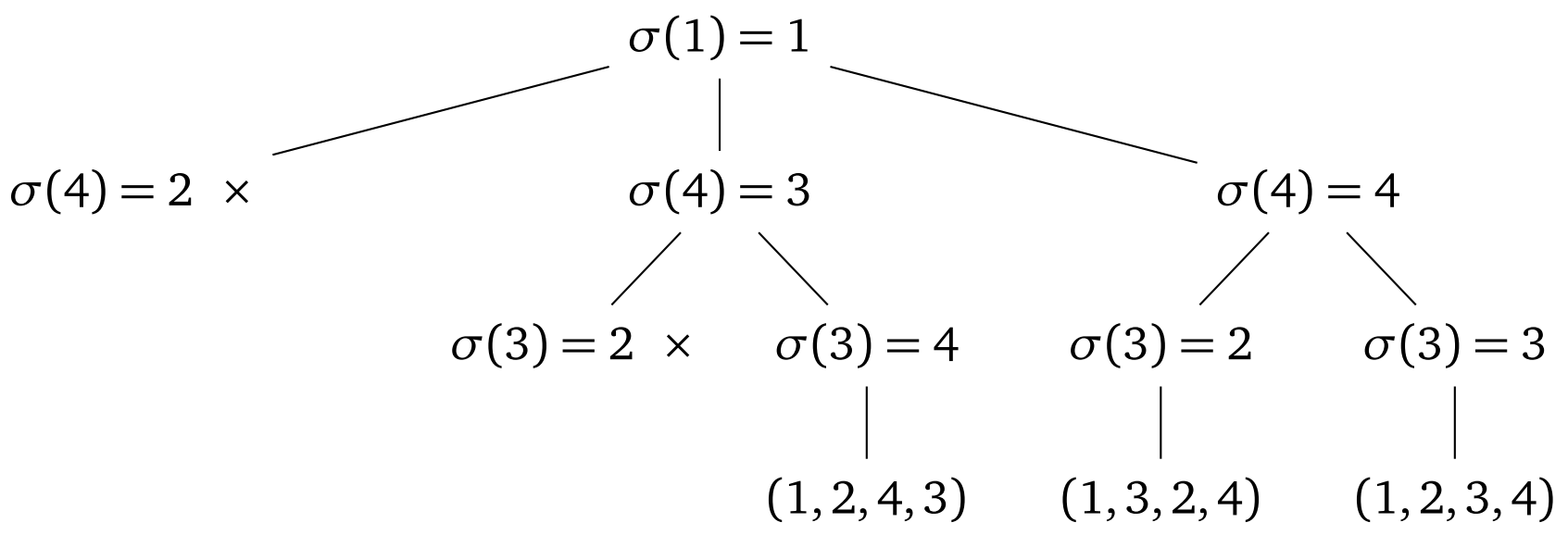

Concretely, the element $\sigma(1) = 1$ is fixed from the beginning. The remaining positions are filled from right to left: at each stage, a candidate element $R$ is appended to the current partial permutation only if the updated running sum (\verb|current_sum += (R-(p+1))|) remains non-negative. If it does not, that branch of the search tree is removed entirely, and no permutations extending that partial assignment are ever constructed (see Fig.~\ref{fig:pruning}). When the remaining pool of numbers $\{2,3,\cdots,N\}$ is exhausted, the completed suffix is prepended with $1$ to yield a valid permutation.

The two algorithms are equivalent in output, as they both implement the same pair of conditions. However, the approach in version 2 prunes the search tree early and avoids the $O(N!)$ memory cost of full enumeration, unlike the approach in version 1. For the purposes of obtaining the set \(\Phi_p\) (and eventually computing $\consta(p)$), version 2 of the algorithm is strictly preferable, and is what makes the computation possible up to $p = 5$ on modest hardware.

\begin{ex}\label{ex:p1p2}
We illustrate the filtering procedure for $p=1$ and $p=2$.

\textbf{Case $p=1$, $N=2$.} The group $\BBS_2$ consists
of two permutations, and thus the alternating sum \eqref{eq:reduced} reads $\sum\nolimits_{\sigma \in \BBS_2} (-1)^\sigma\, \bigl(w_{\sigma(1)}\partial_x\bigr)\bigl(
w_{\sigma(2)}\bigr).$ The two terms are as follows. For $\sigma = (1,2)$, we obtain the term $1 \cdot \partial_x(x)$, which does contribute. For $\sigma = (2,1)$, we obtain the term $x \cdot \partial_x(1)$, which does not contribute because \(\partial_x(1)=0\). This is also consistent with Lemma~\ref{lem:sigma1}, since $\sigma(1) = 2 \neq 1$, this permutation is discarded. Thus $\Phi_1 = \{(1,2)\} \subset \BBS_2,\ \ |\Phi_1| = 1 = \tfrac{1}{2}\,|\BBS_2|.$

\noindent\textbf{Case $p=2$, $N=4$.} The group $\BBS_4$ has $4! = 24$ elements. Applying Lemma~\ref{lem:sigma1} immediately reduces it to the $3! = 6$ permutations with $\sigma(1)=1$. Of these, the non-negativity condition on the exponents $T_k$ eliminates three more, leaving the valid set
\[
\Phi_2 = \{(1,2,3,4),\ (1,2,4,3),\ (1,3,2,4)\} \subset \BBS_4,\ \ |\Phi_2| = 3 = \tfrac{1}{8}\,|\BBS_4|.
\]
The three discarded permutations with $\sigma(1)=1$ are $(1,3,4,2)$, $(1,4,2,3)$, and $(1,4,3,2)$, each of them produces a negative exponent at some intermediate stage and therefore contributes a zero term to the sum.
\end{ex}
\begin{rem}
The permutations that are `close to' the identity permutation $(1, 2, \dots, N)$ contribute to the sum, while permutations that are `close to' the reversed order permutation $(N, N-1, \dots, 2, 1)$ do not contribute to the sum. Indeed, permutations near the identity present the weight functions in increasing order of their degree, ensuring that each intermediate exponent remains sufficiently large to survive $p$ differentiations. On the other hand, permutations near the reversed order present lower-degree weights in the innermost operators, which has a higher chance of driving the exponent below zero.
In Section~\ref{sec:deviation} we measure this closeness to the identity permutation by the deviation vector.
\end{rem}

\begin{table}[H]
\caption{For each $p$, the size of the contributing set
$\Phi_p \subsetneq \BBS_{N}$ with $N = 2p$, its ratio to $(N-2)!$, and the difference between the numbers of even and odd
permutations in $\Phi_p$. }\label{tab:phip}
\centering
\small
\renewcommand{\arraystretch}{1.25}
\begin{tabular}{crrc}
\hline\hline
$p$ & $|\Phi_p|$ & $|\Phi_p|/(N-2)!$ & $e_p - o_p$ \\
\hline
$1$ & $1$ & $1$ & $+1$ \\
$2$ & $3$ & $1.5$ & $-1$ \\
$3$ & $35$ & $\approx 1.458$ & $+1$ \\
$4$ & $1001$ & $\approx 1.390$ & $-1$ \\
$5$ & $53109$ & $\approx 1.317$ & $+1$ \\
$6$ & $4605271$ & $\approx 1.269$ & $-1$ \\
$7$ & $589809987$ & $\approx 1.231$ & $+1$ \\
$8$ & $104899483845$ & $\approx 1.203$ & $-1$ \\
$9$ & $24698838710457$ & $\approx 1.180$ & $+1$ \\
$10$ & $7444522779300351$ & $\approx 1.163$ & $-1$ \\
$11$ & $2792467448952670671$ & $\approx 1.148$ & $+1$ \\
$12$ & \seqsplit{1276369340371154144337} & $\approx 1.136$ & $-1$ \\
$13$ & \seqsplit{698008560075731437652425} & $\approx 1.125$ & $+1$ \\
\hline\hline
\end{tabular}
\end{table}

\begin{rem}
The sequence $|\Phi_p| = 1, 3, 35, 1001, 53109, 4605271, \ldots$ for
$p = 1, 2, 3, 4, 5, 6$ coincides with the terms of odd index $N' = 2p-1$ of
the sequence \seqnum{A147681} in the {\it On-Line Encyclopedia of Integer
Sequences} \cite{oeisA147681}, which counts \emph{late-growing}
permutations\footnote{\label{FootDefLateGrow}\postdisplaypenalty=10000 
By definition, a permutation $\sigma \in \BBS_n$ is \emph{late\/-\/growing} if its
partial sums never exceed the corresponding partial sums of the average
permutation, that is,
\[
\frac{1}{k}\sum\nolimits_{i=1}^{k} \sigma(i) \;\leqslant\; \frac{n+1}{2}
\qquad \text{for all $k \leqslant n$.}
\]
The bound on the right\/-\/hand side is the arithmetic mean of the numbers
$1$,\ $\ldots$,\ $n$; thus the ``average permutation'' has partial sums
$k(n+1)/2$ as $k\geqslant 1$ grows, so dividing by $k$ gives the stated inequality.}
of $\{1, \ldots, n\}$.

Some statistics about the structure of the sets $\Phi_p$ is summarized in Table~\ref{tab:phip}.
We discover that our strategy to restrict the summation to
the subset $\Phi_p \subsetneq \BBS_N$ is effective: the size $|\Phi_p|$ is comparable with $(N-2)!$, whereas $|\BBS_N|=N!$; the ratios $|\Phi_p|/(N-2)!$ decrease monotonically for $2\leqslant p\leqslant 13$.
\end{rem}

\begin{rem}\label{rem:lategrow}
Let $\sigma\in\BBS_N$ be such that $\sigma(1)=1$, and put $\tau(i):=\sigma(i+1)-1$ for $i=1,\ldots,N-1$. The conditions $T_k\geqslant 0$ for~$\sigma$ are the inequalities $\sum_{i=1}^{m}\tau(i)\leqslant m\,p$, where $m=1,\ldots,N-2$, of late growth for the permutation~$\tau$ of $\{1,\ldots,2p-1\}$ (see footnote~\ref{FootDefLateGrow}); the permutations $\sigma$ and~$\tau$ have equal parity. The numbers $e_p$ and~$o_p$ of even and odd permutations in~$\Phi_p$ differ by one, $e_p-o_p=(-1)^{p+1}$, see the last column of Table~\ref{tab:phip}; this is proved in~\cite{OddEven}.
\end{rem}

\section{A closed-form expression for \texorpdfstring{$\consta(p)$}{const(p)}}
\label{sec:closedform}

\noindent In the previous section, we reduced the alternating sum to the summation over the filtered set $\Phi_p$, that is,
\begin{multline*}
\sum_{\sigma \in \BBS_N} (-1)^\sigma\, \bigl(w_{\sigma(1)}\partial_x^p \circ
\cdots \circ w_{\sigma(N-1)}\partial_x^p\bigr)\bigl(w_{\sigma(N)}\bigr)
\\ {}=
\sum_{\sigma \in \Phi_p \subset \BBS_N} (-1)^\sigma\,
\bigl(w_{\sigma(1)}\partial_x^p \circ
\cdots \circ w_{\sigma(N-1)}\partial_x^p\bigr)\bigl(w_{\sigma(N)}\bigr).
\end{multline*}
\begin{rem}\label{rem:monomialchoice}
The equality above holds for the special choice of monomial weights $w_1 = 1$,\ $\ldots$,\ $w_N = x^{N-1}$; the set \(\Phi_p\) was constructed in reference to this choice.
\end{rem}

\noindent We now derive a closed-form expression for each term in this sum, yielding an explicit formula for $\consta(p)$.
\begin{notation}
For any integer \(k\) from the range \(1\leqslant k\leqslant N-1\), and a permutation \(\sigma\in\Phi_p\), we, by definition, put
\begin{equation}\label{eq:Ek}
E_k \;:=\; \sum_{j=0}^{k-1}\bigl(\sigma(N-j)-1\bigr) - (k-1)p,
\end{equation}
so that its falling factorial is
\[
E_k^{\,\underline{p}}= E_k!/(E_k - p)!.
\]
\end{notation}
\begin{rem}\label{rem:TE}
For $k=1,\ldots,N-1$, the partial sums $T_k$ from \S\ref{sec:algorithm} and the numbers $E_k$ in~\eqref{eq:Ek} differ by~$p$, namely $T_k=E_k-p$. Indeed, $E_k$ is the exponent of~$x$ before the $k$th differentiation $\partial_x^p$, counted from the right, and $T_k$ is the exponent after it. Hence the condition $T_k\geqslant 0$ is $E_k\geqslant p$. For $k=N-1$ we have $T_{N-1}=1-\sigma(1)$ because $(N-1)\,p=\tfrac{1}{2}N(N-1)$, so that the condition $T_{N-1}\geqslant 0$ is the equality $\sigma(1)=1$ of Lemma~\ref{lem:sigma1}.
\end{rem}
\begin{claim}\label{claim:closedform}
\begin{equation}\label{eq:closedform}
\consta(p) \cdot \prod_{k=0}^{N-1} k!
\;=\;
\sum_{\sigma \in \Phi_p} (-1)^\sigma \prod_{k=1}^{N-1}
E_k^{\,\underline{p}}
\end{equation}
\end{claim}

\begin{proof}
For our chosen set of monomial weights, each surviving term in the l.-h.s.\ of \eqref{eq:reduced} is
\begin{multline*}
\left(1\cdot\partial_x^p\circ w_{\sigma(2)}\cdot\partial_x^p\circ\dots\circ w_{\sigma(N-1)}\cdot\partial_x^p\right)(w_{\sigma(N)})\\=\left(1\cdot\partial_x^p\circ x^{\sigma(2)-1}\cdot\partial_x^p\circ\dots\circ x^{\sigma(N-1)-1}\cdot\partial_x^p\right)(x^{\sigma(N)-1}).
\end{multline*}
The r.h.s\ of formula~\eqref{eq:closedform} tracks the scalar coefficient accumulated as the composition of weighted operators acts on the monomial $x^{\sigma(N)-1}$, proceeding from right to left through the $N-1$ applications of $\partial_x^p$. Let us study this right-to-left action of differential operators in more detail.

\noindent\textbf{Rightmost action.} The rightmost operator acts on the weight $w_{\sigma(N)} = x^{\sigma(N)-1}$:
\[
\partial_x^p\!\left(x^{\sigma(N)-1}\right)
= \frac{(\sigma(N)-1)!}{(\sigma(N)-1-p)!}\, x^{\sigma(N)-1-p}
= E_1^{\,\underline{p}}\cdot x^{\sigma(N)-1-p},
\]
where $E_1 = \sigma(N)-1$. The exponent of $x$ is now $\sigma(N)-1-p$.

\noindent\textbf{Second action.} Pre-multiplying by $w_{\sigma(N-1)} = x^{\sigma(N-1)-1}$ shifts the exponent to
\[
(\sigma(N-1)-1) + (\sigma(N)-1-p)
= \sigma(N-1)+\sigma(N)-2-p = E_2,
\]
and a further application of $\partial_x^p$ yields the factor $E_2!/(E_2-p)! = E_2^{\,\underline{p}}$.
\textbf{General step.} After $k$ actions of $\partial_x^p$, the accumulated exponent of $x$ is $E_k - p$. Pre-multiplying by $w_{\sigma(N-k)} = x^{\sigma(N-k)-1}$ raises the exponent to $E_{k+1}$, and the $(k+1)$-th application of $\partial_x^p$ contributes the factor $E_{k+1}^{\,\underline{p}}$. Iterating through all the $N-1$ actions, the total scalar coefficient accumulated for the permutation $\sigma\in\Phi_p$ is
\[
c(\sigma,p) \;=\; \prod_{k=1}^{N-1} E_k^{\,\underline{p}}.
\]
Summing over all the contributing permutations with their signs and equating it to the product of $\consta(p)$ and the Wronskian determinant $\Wronsk(1,x,\dots,x^{N-1}) = \prod_{k=0}^{N-1}k!$ from Lemma~\ref{lem:wronskian} yields equality \eqref{eq:closedform}.
\end{proof}

\begin{ex}
\label{ex:peq2}
We now illustrate Claim~\ref{claim:closedform} for $p=2$, $N=4$. The contributing permutations, as identified in Example~\ref{ex:p1p2}, are $\Phi_2 = \{(1,2,3,4),\, (1,2,4,3),\, (1,3,2,4)\}$. For each of them we compute the sign and the coefficient $c(\sigma, 2) = \prod_{k=1}^{3} E_k^{\,\underline{2}}$:
\[
\renewcommand{\arraystretch}{1.2}
\begin{array}{cclll}
\sigma & (-1)^\sigma & E_1 & E_2 & E_3 \qquad\ c(\sigma,2)\\
\hline
(1,2,3,4) & +1 & 3, & 3+2-2=3, & 3+1-2=2; \quad 3^{\underline{2}}\cdot 3^{\underline{2}}\cdot 2^{\underline{2}} = 3!\cdot 3!\cdot 2!\\
(1,2,4,3) & -1 & 2, & 2+3-2=3, & 3+1-2=2; \quad 2^{\underline{2}}\cdot 3^{\underline{2}}\cdot 2^{\underline{2}} = 3!\cdot 2!\cdot 2!\\
(1,3,2,4) & -1 & 3, & 3+1-2=2, & 2+2-2=2; \quad 3^{\underline{2}}\cdot 2^{\underline{2}}\cdot 2^{\underline{2}} = 3!\cdot 2!\cdot 2!
\end{array}
\]
Assembling the signed sum and dividing by the Wronskian determinant $\Wronsk(1,x,x^2,x^3)$ $= 3!\cdot 2!\cdot 1!$ from Lemma~\ref{lem:wronskian}, we obtain that
\[
\consta(2) = \frac{3!\cdot 3!\cdot 2! - 3!\cdot 2!\cdot 2! - 3!\cdot 2!\cdot 2!}{3!\cdot 2!\cdot 1!} = 3! - 2! - 2! = 2,
\]
in agreement with the value \(\consta(2)=2\) recorded in~\eqref{eq:smallp}.
\end{ex}

\section{An intermediate acceleration: the deviation-vector formulation}
\label{sec:deviation}

\noindent The filtering procedure of Section~\ref{sec:algorithm} associates to each permutation a binary outcome: either the permutation belongs to \(\Phi_p\) or it does not. The following deviation-vector formulation refines this viewpoint by assigning to each contributing permutation a vector of nonnegative integers; in doing so we measure its distance from the extremal case given by the identity permutation. In this sense, it builds on our qualitative condition by also measuring how strongly a contributing permutation \(\sigma\in\Phi_p\) satisfies the late-growing condition.

The motivation to introduce this statistic on late-growing permutations and provide the reformulation below was to develop an improvement of the original algorithm. Although it yields a significant speedup and was sufficient to obtain \(\consta(7)\), it was later superseded by a fully recursive method with memoization~\cite{Method}. We nevertheless present it because we now obtain a mathematically whole idea of the role of late-growing permutations in our computations.

\begin{define}\label{def:master}
Let $\sigma_{\mathrm{id}} = (1,2,\ldots,N)$ denote the identity permutation, which always
belongs to $\Phi_p$. Substituting $\sigma_{\mathrm{id}}$ into~\eqref{eq:Ek}, the $k$-th
component of the resulting vector is
\[
v^{\mathrm{m}}_k
\;:=\; E_k\!\left(\sigma_{\mathrm{id}}\right)
\;=\; \sum_{j=0}^{k-1}(N - 1 - j) \;-\; (k-1)p
\;=\; k(N-1) - \tfrac{k(k-1)}{2} - (k-1)p.
\]
We call $\mathbf{v}^{\mathrm{m}} = (v^{\mathrm{m}}_1, \ldots, v^{\mathrm{m}}_{N-1})$
the \emph{master vector}; it is determined solely by $p$.
\end{define}

\begin{lemma}\label{lem:masterclosed}
For $N=2p$ and $k=1,\dots,N-1$ the `master vector' $v^{\mathrm m}_k$ equals
\[
v^{\mathrm m}_k \;=\; (k+1)(2p-k)/2.
\]
Consequently, over $k\in\{1,\dots,N-1\}$, $v^{\mathrm m}_k$ attains its maximum at
$k=p-1$ and at $k=p$, where $\max_{1\leqslant k\leqslant N-1} v^{\mathrm m}_k = p(p+1)/2$.
\end{lemma}

\begin{proof}
Substituting $N=2p$ into $v^{\mathrm m}_k = k(N-1)-\tfrac{k(k-1)}{2}-(k-1)p$ and
simplifying gives
\[
v^{\mathrm m}_k = k(2p-1)-\tfrac{1}{2}k(k-1)-(k-1)p
= p(k+1) - \tfrac{1}{2}k(k+1) = \tfrac{1}{2}(k+1)(2p-k).
\]
The quadratic polynomial $(k+1)(2p-k)$ in~$k$ attains its maximum at $k=p-\tfrac12$, hence over the integers at $k=p-1$ and at $k=p$, where it equals $p(p+1)$.
\end{proof}

\begin{define}\label{def:deviation}
For any $\sigma\in\Phi_p$, define the \emph{deviation vector}
$\boldsymbol{\delta}(\sigma) = (\delta_1(\sigma),\ldots,\delta_{N-1}(\sigma))$ by
\[
\delta_k(\sigma)
\;:=\; v^{\mathrm{m}}_k - E_k(\sigma)
\;=\; \sum_{j=0}^{k-1}\bigl[(N - j) - \sigma(N-j)\bigr].
\]
\end{define}

\begin{rem}\label{rem:deltanonneg}
Each $\delta_k(\sigma) \geqslant 0$ tracks how much the right-to-left running sum of $\sigma$
falls short of the identity; it is nonnegative because the $k$ rightmost entries of~$\sigma$ are distinct elements of $\{1,\ldots,N\}$, whence their sum does not exceed $N+(N-1)+\cdots+(N-k+1)$. In particular, $\delta_k = 0$ for all $k$ if and only if
$\sigma = \sigma_{\mathrm{id}}$.
\end{rem}

\begin{define}\label{def:oneline}
Let $\lambda\colon\BBZ[\BBS_N]\to\BBZ^N$ be the $\BBZ$-linear map which takes every permutation $\sigma\in\BBS_N$ to its one\/-\/line notation, $\lambda(\sigma):=(\sigma(1),\ldots,\sigma(N))$. Let $\rho\colon\BBZ^N\to\BBZ^{N-1}$ be the $\BBZ$-linear map of partial sums from the right: for $u=(u_1,\ldots,u_N)$ and $k=1,\ldots,N-1$,
\[
\rho(u)_k := \sum_{j=0}^{k-1} u_{N-j}.
\]
\end{define}

\begin{lemma}\label{lem:devgroupalg}
For every $\sigma\in\Phi_p$, the deviation vector is the image of the element $\id-\sigma$ of the group algebra $\BBZ[\BBS_N]$:
\begin{equation}\label{eq:devgroupalg}
\boldsymbol{\delta}(\sigma) = \rho\bigl(\lambda(\id-\sigma)\bigr).
\end{equation}
\end{lemma}

\begin{proof}
By linearity, $\lambda(\id-\sigma)=\lambda(\id)-\lambda(\sigma)=\bigl(1-\sigma(1),\ldots,N-\sigma(N)\bigr)$ is the termwise difference of the one\/-\/line notations. By definition, its $k$th partial sum from the right is $\delta_k(\sigma)$.
\end{proof}

\begin{define}\label{def:F}
The \emph{product function} is the map $F\colon \mathbb{Z}_{\geqslant p}^{N-1} \to \mathbb{Z}_{>0}$ given by
\[
F(\mathbf{v})
\;:=\; \prod_{k=1}^{N-1} \frac{v_k!}{(v_k - p)!}
\;=\; \prod_{k=1}^{N-1} v_k^{\,\underline{p}}.
\]
\end{define}

\begin{prop}\label{prop:deviation}
Formula~\eqref{eq:closedform} can be written as
\begin{equation}\label{eq:deviation}
\consta(p) \cdot \prod_{k=0}^{N-1} k!
\;=\;
\sum_{\sigma\in\Phi_p} (-1)^\sigma\,
F\!\left(\mathbf{v}^{\mathrm{m}} - \boldsymbol{\delta}(\sigma)\right).
\end{equation}
\end{prop}

\begin{proof}
By Definition~\ref{def:deviation}, $E_k(\sigma)=v^{\mathrm{m}}_k-\delta_k(\sigma)$ for every $\sigma\in\Phi_p$ and every~$k$.
\end{proof}

\begin{rem}\label{rem:domainF}
The filtering condition $T_k \geqslant 0$ (which defines $\Phi_p$) is precisely the condition
that $E_k(\sigma) \geqslant p$ for all $k$ (see Remark~\ref{rem:TE}), i.e.\ $\delta_k(\sigma) \leqslant v^{\mathrm{m}}_k - p$.
The filter thus guarantees that all arguments of $F$ remain in the domain $\mathbb{Z}_{\geqslant p}$.
\end{rem}

After this reformulation, the entire computation of $\consta(p)$ reduces to: generate $\Phi_p$ incrementally, read off $\boldsymbol{\delta}(\sigma)$ for free during generation, evaluate $F$,
and accumulate the signed sum (see Algorithm~\ref{alg:dev}).

\begin{algorithm}[htb]
\caption{The deviation-vector evaluation of the right-hand side of~\eqref{eq:deviation}.}\label{alg:dev}
\begin{algorithmic}[1]
\Require $p\geqslant 1$; put $N:=2p$, $v^{\mathrm m}_k:=\tfrac12(k+1)(2p-k)$ and $W:=\prod_{k=0}^{N-1}k!$; tabulate $E^{\,\underline{p}}$ for $p\leqslant E\leqslant\tfrac12 p(p+1)$
\State $S\gets 0$
\Procedure{Extend}{$k,\ U,\ \boldsymbol{\delta},\ \varepsilon$}
  \Comment{$\sigma(N-k+1),\ldots,\sigma(N)$ are placed}
  \Statex \Comment{$U\subseteq\{2,\ldots,N\}$: the entries not yet placed;\ \ $\boldsymbol{\delta}=(\delta_1,\ldots,\delta_k)$, $\delta_0:=0$;\ \ $\varepsilon=\pm1$}
  \If{$k=N-1$}
    \State $S\gets S+\varepsilon\cdot\prod_{j=1}^{N-1}\bigl(v^{\mathrm m}_j-\delta_j\bigr)^{\underline{p}}$
      \Comment{$=(-1)^\sigma F(\mathbf{v}^{\mathrm m}-\boldsymbol{\delta}(\sigma))$, by table look-ups}
    \State \Return
  \EndIf
  \For{$R\in U$}
    \State $d\gets\delta_k+(N-k)-R$
      \Comment{put $\sigma(N-k):=R$; this is $\delta_{k+1}$, see~\eqref{eq:devgroupalg}}
    \If{$d\leqslant v^{\mathrm m}_{k+1}-p$}
      \Comment{otherwise the branch is pruned}
      \State \Call{Extend}{$k+1,\ U\setminus\{R\},\ (\boldsymbol{\delta},d),\ \varepsilon\cdot(-1)^{\#\{R'\in\{2,\ldots,N\}\setminus U\,:\,R'<R\}}$}
    \EndIf
  \EndFor
\EndProcedure
\State \Call{Extend}{$0,\ \{2,\ldots,N\},\ (\,),\ +1$}
\State \Return $S/W$
  \Comment{$=\consta(p)$ by~\eqref{eq:deviation}}
\end{algorithmic}
\end{algorithm}
\begin{ex}
Consider $p=3$, $N=6$; we have that $\mathbf{v}^{\mathrm{m}} = (5,6,6,5,3)$. This vector corresponds exactly to the factorial arguments appearing in the leading (identity) term of the expansion. Naturally, every admissible permutation contributes a term obtained by lowering certain entries of the master vector. For instance,
\[
\renewcommand{\arraystretch}{1.6}
\begin{array}{cccl}
\sigma & \boldsymbol{\delta}(\sigma) & \mathbf{v}^{\mathrm{m}}-\boldsymbol{\delta}(\sigma) & \text{contribution}\\
\hline
(1,2,3,4,5,6) & (0,0,0,0,0) & (5,6,6,5,3) & +\,\frac{5!}{2!}\cdot\frac{6!}{3!}\cdot\frac{6!}{3!}\cdot\frac{5!}{2!}\cdot\frac{3!}{0!}\\
(1,2,3,4,6,5) & (1,0,0,0,0) & (4,6,6,5,3) & -\,\frac{4!}{1!}\cdot\frac{6!}{3!}\cdot\frac{6!}{3!}\cdot\frac{5!}{2!}\cdot\frac{3!}{0!}\\
(1,2,3,5,4,6) & (0,1,0,0,0) & (5,5,6,5,3) & -\,\frac{5!}{2!}\cdot\frac{5!}{2!}\cdot\frac{6!}{3!}\cdot\frac{5!}{2!}\cdot\frac{3!}{0!}
\end{array}
\]
\end{ex}
\begin{rem}
For small values of $p$, the master vectors are: $\mathbf v^{\mathrm m}=(1)$ for $p=1$; $\mathbf v^{\mathrm m}=(3,3,2)$ for $p=2$; $\mathbf v^{\mathrm m}=(5,6,6,5,3)$ for $p=3$; $\mathbf v^{\mathrm m}=(7,9,10,10,9,7,4)$ for $p=4$.
\end{rem}

\medskip
\noindent\textbf{Timings.}
The run times of the generators ver.~1 and ver.~2 from \S\ref{sec:algorithm} and of Algorithm~\ref{alg:dev} are compared in Table~\ref{tab:timings}.
\begin{table}[htb]
\caption{Run times, in seconds, of three evaluations of the sum $\consta(p)\cdot\prod_{k=0}^{N-1}k!$\,: the generator ver.~1 and the generator ver.~2, each followed by formula~\eqref{eq:closedform} for every $\sigma\in\Phi_p$, and Algorithm~\ref{alg:dev}. One thread of an Intel Xeon processor at 2.1\,GHz; the implementations are in \textsf{Python}~3.13 (left) and in \textsf{C++} with 576-bit integers (right). The codes are available at~\cite{Code}. A dash marks a run which was not performed.}\label{tab:timings}
\centering
\small
\setlength{\tabcolsep}{4.5pt}
\renewcommand{\arraystretch}{1.2}
\begin{tabular}{crr|rrr|rrr}
\hline\hline
 & & & \multicolumn{3}{c|}{\textsf{Python} timing} & \multicolumn{3}{c}{\textsf{C++} timing} \\
$p$ & $N!$ & $|\Phi_p|$ & ver.~1 & ver.~2 & Alg.~\ref{alg:dev} & ver.~1 & ver.~2 & Alg.~\ref{alg:dev} \\
\hline
$4$ & $40320$ & $1001$ & $0.00524$ & $0.00325$ & $0.00182$ & $0.00028$ & $0.00016$ & $0.00010$ \\
$5$ & $3628800$ & $53109$ & $0.461$ & $0.258$ & $0.130$ & $0.015$ & $0.00745$ & $0.00499$ \\
$6$ & $479001600$ & $4605271$ & $58.4$ & $28.0$ & $11.6$ & $2.08$ & $0.929$ & $0.519$ \\
$7$ & $87178291200$ & $589809987$ & --- & --- & $1663$ & $371$ & $153$ & $88.5$ \\
\hline\hline
\end{tabular}
\end{table}

\begin{rem}\label{rem:timings}
Every run in Table~\ref{tab:timings} returns the value $\consta(p)$ from Table~\ref{tab:values}. The generator ver.~1 runs through all the $N!$ permutations, whereas ver.~2 and Algorithm~\ref{alg:dev} reach only the $|\Phi_p|$ permutations which contribute; at $p=7$ these are $N!/|\Phi_p|\approx 148$ times fewer. In \textsf{C++} at $p=7$, Algorithm~\ref{alg:dev} is $4.2$ times faster than ver.~1 and $1.7$ times faster than ver.~2; in \textsf{Python} at $p=6$, these factors are $5.0$ and $2.4$, respectively.
\end{rem}

\section{New values of the new sequence, and its growth rate}\label{sec:results}

\begin{table}[htb]
\caption{Values of $\consta(p)$: exact for $p\leqslant 7$, approximate for $8\leqslant p\leqslant 20$ (for $p\leqslant 18$ the latter are taken from~\cite{Method}).}\label{tab:values}
\centering
\small
\renewcommand{\arraystretch}{1.2}
\begin{tabular}{rl@{\qquad\qquad}rl}
\hline\hline
$p$ & $\consta(p)$ & $p$ & $\consta(p)$ \\
\hline
$1$ & $1$ & $11$ & $\approx 1.512\cdot 10^{126}$ \\
$2$ & $2$ & $12$ & $\approx 5.498\cdot 10^{159}$ \\
$3$ & $90$ & $13$ & $\approx 8.396\cdot 10^{197}$ \\
$4$ & $586656$ & $14$ & $\approx 7.398\cdot 10^{240}$ \\
$5$ & $1915103977500$ & $15$ & $\approx 5.046\cdot 10^{288}$ \\
$6$ & $7886133184567796056800$ & $16$ & $\approx 3.503\cdot 10^{341}$ \\
$7$ & $85873408332103907284746052081828368$ & $17$ & $\approx 3.198\cdot 10^{399}$ \\
$8$ & $\approx 4.594\cdot 10^{51}$ & $18$ & $\approx 4.881\cdot 10^{462}$ \\
$9$ & $\approx 2.060\cdot 10^{72}$ & $19$ & $\approx 1.563\cdot 10^{531}$ \\
$10$ & $\approx 1.237\cdot 10^{97}$ & $20$ & $\approx 1.301\cdot 10^{605}$ \\
\hline\hline
\end{tabular}
\end{table}

\begin{claim}\label{ClaimPositiveGrows}
(\textit{i}), in~\cite[\S2]{Dzhuma2014}: The universal constants $\consta(p)$ in equality~\eqref{EqMain} are positive: $\consta(p)>0$ for all $p\in\BBN_{\geqslant 1}$.\\[0.5pt]
(\textit{ii}), by O.\,Zaboronski${}^{\text{\protect\ref{FootZabor}}}$: The integer sequence $\bigl\{ \consta(p) \bigr\}_{p\in\BBN}$ grows monotonically: $1 \leqslant \consta(p) < \consta(q)$ if $1\leqslant p < q$.
\end{claim}

\subsection{Prime factorisations of known \texorpdfstring{$\consta(p)$}{const(p)}}\label{sec:factorisations}

\begin{notation}\label{not:padic}
For a prime $p$ and a nonzero rational number $r$, let $v_p(r)\in\mathbb{Z}$ denote its $p$-adic valuation: writing $r=p^{e}\cdot\frac{a}{b}$ with $a,b$ integers coprime
to $p$, we set $v_p(r):=e$.
\end{notation}

\begin{table}[htb]
\caption{The $q$-adic valuations $v_q(\consta(p))$ for primes
$q \leqslant 19$ and $2 \leqslant p \leqslant 20$ (for $8\leqslant p\leqslant 18$ they are obtained from the exact values in~\cite{Method}). The entries $v_p(\consta(p))$
for prime $p$ are in boldface.}
\label{tab:valuations}
\centering
\footnotesize
\setlength{\tabcolsep}{3.1pt}
\renewcommand{\arraystretch}{1.2}
\begin{tabular}{r*{19}{c}}
\hline\hline
$q \,\backslash\, p$
& 2 & 3 & 4 & 5 & 6 & 7 & 8 & 9 & 10 & 11 & 12 & 13 & 14 & 15 & 16 & 17 & 18 & 19 & 20 \\
\hline
2 & $\mathbf{1}$ & 1 & 5 & 2 & 5 & 4 & 17 & 11 & 16 & 7 & 18 & 13 & 17 & 11 & 49 & 33 & 32 & 16 & 38\\
3 & $\cdot$ & $\mathbf{2}$ & 3 & 1 & 6 & 2 & 2 & 14 & 6 & 5 & 11 & 8 & 3 & 12 & 3 & 3 & 31 & 15 & 9\\
5 & $\cdot$ & 1 & $\cdot$ & $\mathbf{4}$ & 2 & $\cdot$ & $\cdot$ & $\cdot$ & 8 & 6 & 5 & 2 & 2 & 12 & 6 & 1 & 2 & 1 & 16\\
7 & $\cdot$ & $\cdot$ & 1 & 1 & $\cdot$ & $\mathbf{6}$ & 5 & 2 & $\cdot$ & $\cdot$ & $\cdot$ & $\cdot$ & 13 & 8 & 4 & 2 & $\cdot$ & $\cdot$ & $\cdot$\\
11 & $\cdot$ & $\cdot$ & $\cdot$ & $\cdot$ & 1 & $\cdot$ & $\cdot$ & 2 & $\cdot$ & $\mathbf{10}$ & 8 & 6 & 4 & 3 & 1 & $\cdot$ & $\cdot$ & $\cdot$ & $\cdot$\\
13 & $\cdot$ & $\cdot$ & $\cdot$ & $\cdot$ & $\cdot$ & 1 & 1 & 1 & $\cdot$ & 1 & $\cdot$ & $\mathbf{12}$ & 11 & 8 & 6 & 5 & 2 & 1 & $\cdot$\\
17 & $\cdot$ & $\cdot$ & $\cdot$ & $\cdot$ & $\cdot$ & $\cdot$ & $\cdot$ & 1 & $\cdot$ & $\cdot$ & $\cdot$ & 1 & $\cdot$ & 1 & $\cdot$ & $\mathbf{16}$ & 15 & 12 & 10\\
19 & $\cdot$ & $\cdot$ & $\cdot$ & $\cdot$ & $\cdot$ & $\cdot$ & $\cdot$ & $\cdot$ & 1 & 1 & $\cdot$ & $\cdot$ & $\cdot$ & $\cdot$ & $\cdot$ & 1 & $\cdot$ & $\mathbf{18}$ & 17\\
\hline\hline
\end{tabular}
\end{table}

\subsection{Growth rate}
\label{sec:growth}

\begin{claim}[{O.\,Zaboronski${}^{\text{\protect\ref{FootZabor}}}$}]\label{ClaimAsympt}
As $p\to +\infty$, the asymptotic expansion is
\begin{equation}\label{EqAsymptZabor}
\log\consta(p) = 2p^2\log p - b p^2 + \overline{o}(p^2),
\end{equation}
where $b\geqslant 2.6744$.
\end{claim}

\section{A proven lower bound on \texorpdfstring{$v_p(\consta(p))$}{v\_p(const(p))}}
\label{sec:lowerbound}\label{sec:padic}

\noindent In this section we will prove that for prime \(p\), \(v_p(\consta(p))\geqslant p-1\), which gives a lower bound on $v_p(\consta(p))$.

We recall the ultrametric inequality for the $p$-adic valuation (Notation~\ref{not:padic}),
\[
v_p(x+y)\;\geqslant\;\min\bigl(v_p(x),v_p(y)\bigr),
\]
valid for all rationals $x,y$; it extends inductively to any
finite sum\footnote{
\(
v_p\!\left(\sum_{i=1}^{n} x_i\right)\geqslant \min_{1\leqslant i\leqslant n} v_p(x_i)
\)
for all finite collections \(\{x_i\}\). Indeed, the case \(n=2\) is the ultrametric inequality. If the result holds for \(n-1\) terms, then
\[
v_p\!\left(\sum_{i=1}^{n}x_i\right)=
v_p\!\left(\left(\sum_{i=1}^{n-1}x_i\right)+x_n\right)
\geqslant
\min\!\left(
v_p\!\left(\sum_{i=1}^{n-1}x_i\right),
v_p(x_n)
\right)
\geqslant
\min_{1\leqslant i\leqslant n} v_p(x_i),
\]
which completes the induction.
}.

\begin{lemma}\label{lem:Erange}
For every $\sigma\in\Phi_p$ and every $k=1,\dots,N-1$,
\[
p \;\leqslant\; E_k(\sigma) \;\leqslant\; \tfrac{1}{2}p(p+1) \;<\; p^2.
\]
\end{lemma}

\begin{proof}
The lower bound is precisely the defining condition of $\Phi_p$ (\S\ref{sec:algorithm} and Remark~\ref{rem:TE}); equivalently,
$\delta_k(\sigma)\leqslant v^{\mathrm m}_k-p$, as noted following~\eqref{eq:deviation}. The upper bound follows from $\delta_k(\sigma)\geqslant 0$ together with Lemma~\ref{lem:masterclosed}:
\[
E_k(\sigma)=v^{\mathrm m}_k-\delta_k(\sigma)\leqslant
v^{\mathrm m}_k\leqslant\tfrac{1}{2}p(p+1).
\] The strict inequality $\tfrac{1}{2}p(p+1)<p^2$ is equivalent to $p>1$, which holds for every prime $p$.
\end{proof}

\begin{lemma}\label{lem:fallingvaluation}
Let $p$ be a prime and let $m$ be an integer with $p\leqslant m<p^2$. Then\linebreak
$v_p\bigl(m!/(m-p)!\bigr) = 1$.
\end{lemma}

\begin{proof}
We recall Legendre's formula, which states that $v_p(n!)=\sum_{i\geqslant1}\lfloor n/p^i\rfloor$ for $n\geqslant0$. Both $m$
and $m-p$ lie in $[0,p^2)$, so every term with $i\geqslant2$ vanishes, giving
$v_p(m!)=\lfloor m/p\rfloor$ and $v_p((m-p)!)=\lfloor (m-p)/p\rfloor=\lfloor m/p\rfloor-1$.
Subtracting the two expressions gives
\[
v_p\bigl(m!/(m-p)!\bigr)=v_p(m!)-v_p((m-p)!)=1,
\]
as claimed.
\end{proof}

\begin{lemma}\label{lem:cvaluation}
For every prime $p$ and every $\sigma\in\Phi_p$,
\[
v_p\bigl(c(\sigma,p)\bigr) \;=\; N-1 \;=\; 2p-1,
\qquad\text{where } c(\sigma,p)=\prod_{k=1}^{N-1}E_k^{\,\underline{p}}.
\]
\end{lemma}

\begin{proof}
By Lemma~\ref{lem:Erange}, each $E_k(\sigma)$ satisfies $p\leqslant E_k(\sigma)<p^2$, so
Lemma~\ref{lem:fallingvaluation} gives $v_p\bigl(E_k(\sigma)^{\,\underline{p}}\bigr)=1$
for every $k=1,\dots,N-1$. Since $v_p$ is additive over products,
\[
v_p(c(\sigma,p)) = \sum_{k=1}^{N-1}v_p\bigl(E_k(\sigma)^{\,\underline{p}}\bigr) = N-1,
\]
as claimed.
\end{proof}

\begin{lemma}\label{lem:Wvaluation}
Let $W:=\prod_{k=0}^{N-1}k!=\Wronsk(1,x,\dots,x^{N-1})$ (Lemma~\ref{lem:wronskian}).
Then $v_p(W) = p$.
\end{lemma}

\begin{proof}
For $0\leqslant k\leqslant p-1$, $v_p(k!)=0$ since $k<p$. For $p\leqslant k\leqslant N-1=2p-1$, we have $k<p^2$
(as $2p-1<p^2\iff(p-1)^2>0$, true for $p\geqslant2$), so Legendre's formula gives
$v_p(k!)=\lfloor k/p\rfloor=1$. Exactly $p$ of the $N=2p$ values of $k$ fall in this
upper range, each contributing $1$ and the rest contributing $0$, so $v_p(W)=p$.
\end{proof}

\begin{theor}\label{thm:lowerbound}
For every prime $p$,
\[
v_p\bigl(\consta(p)\bigr) \;\geqslant\; p-1.
\]
\end{theor}

\begin{proof}
By Claim~\ref{claim:closedform},
\[
\consta(p)\cdot W \;=\; \sum_{\sigma\in\Phi_p}(-1)^\sigma\,c(\sigma,p).
\]
By Lemma~\ref{lem:cvaluation}, every summand on the right satisfies
$v_p\bigl((-1)^\sigma c(\sigma,p)\bigr)=2p-1$, the sign not affecting the valuation.
By the ultrametric inequality,
\[
v_p\!\left(\sum_{\sigma\in\Phi_p}(-1)^\sigma\,c(\sigma,p)\right)\;\geqslant\;
\min_{\sigma\in\Phi_p} \left\{v_p\bigl((-1)^\sigma c(\sigma,p)\bigr)\right\} \;=\; 2p-1.
\]
Thus $v_p(\consta(p)\cdot W)\geqslant 2p-1$, and by Lemma~\ref{lem:Wvaluation},
\[
v_p(\consta(p)) = v_p(\consta(p)\cdot W) - v_p(W) \geqslant (2p-1)-p = p-1,
\]
as claimed.
\end{proof}

One feature of Table~\ref{tab:valuations} is immediately apparent: the prime factors of $p$ appear in $\consta(p)$ with notably high valuation. Particularly, for every prime $p$ in our range,
the entry $v_p(\consta(p))$ equals exactly $p - 1$.

\begin{conjecture}\label{conj:valuation}
For every prime $p$,
\[
v_p\bigl(\consta(p)\bigr) \;=\; p - 1.
\]
\end{conjecture}

\begin{rem}
Theorem~\ref{thm:lowerbound} establishes half of Conjecture~\ref{conj:valuation}
unconditionally for every prime $p$. Writing $\consta(p)=p^{p-1}\cdot S_p$ for the (now
guaranteed) integer $S_p:=\consta(p)/p^{p-1}$, the remaining content of the
conjecture is the single statement $p\nmid S_p$, i.e.\ the absence of any further
cancellation modulo $p$ in the alternating sum of Claim~\ref{claim:closedform} beyond
what is accounted for above.
Note that the proof of Theorem~\ref{thm:lowerbound} uses the permutations $\sigma\in\Phi_p$ only through Lemma~\ref{lem:Erange}, that is, through the inequalities $0\leqslant\delta_k(\sigma)\leqslant v^{\mathrm{m}}_k-p$ for the deviation vector $\boldsymbol{\delta}(\sigma)=\rho\bigl(\lambda(\id-\sigma)\bigr)$.
\end{rem}

\section{Conclusion}
\nopindent%
The subset $\Phi_p\subseteq \BBS_{2p-1} \subsetneq\BBS_{2p}$ of late\/-\/growing permutations (of the $n=2p-1$ objects $\{2,\ldots,2p\}$, with $\sigma(1)\equiv 1$ fixed in the $\BBS_{2p}$-\/description of~$\Phi_p$) is the set of those permutations which do contribute --\,to the count of the Wronskian determinants, see Eq.~\eqref{EqMain}\,-- in the course of alternating $N=2p$ differential operators $1\cdot\dd_x^p$,\ $x\cdot\dd_x^p$,\ $\ldots$,\ $x^{2p-1}\cdot\dd_x^p$ of strict order $p\in\BBN_{\geqslant 1}$ on the line~$\BBR$ with affine coordinate~$x$ (see~\cite{cpWron,Method,OddEven}).
Under the assumption that the operator coefficients are the monomials $w_1=1$,\ $w_2=x$,\ $\ldots$,\ $w_N=x^{N-1}$, the conditions of `late growth' for a permutation $\sigma\in\Phi_p\subseteq \BBS_{2p-1} \subsetneq \BBS_{2p}$, namely 
$\sum_{i=2}^{k+1} \bigl( \sigma(i) - 1 \bigr) \leqslant k\cdot(n+1)/2$ 
for all positive $k\leqslant n$, stems from the ban to have negative powers of~$x$ while the alternated operators $x^{i-1}\cdot\dd_x^p$ act from the right to the left (with the coefficient $w_{\sigma(N)}$ standing rightmost).\footnote{\label{FootAltOrdP}
The alternated sum of operators' associative compositions
is again an operator of strict order~$p$, with the universal constant $\consta(p)$ times the Wronskian $\Wronsk\bigl(w_1,\ldots,w_N\bigr)$ for the new coefficient (see~\cite{cpWron,Method,OddEven}). This is why it suffices to act by the $N-1$ operators $1\cdot\dd_x^p$,\ $x^{\sigma(2)-1}\cdot\dd_x^p$,\ $\ldots$,\ $x^{\sigma(N-1)-1}\cdot\dd_x^p$ on the rightmost coefficient~$w_{\sigma(N)}$.}
By alternating over these much smaller (than $\BBS_{2p}$) subsets~$\Phi_p$, we develop an algorithm to count the Wronskians $\Wronsk\bigl(1,x,\ldots,x^{2p-1}\bigr)
= \prod_{k=0}^{2p-1} k!$ as those appear, with the universal constant $\consta(p)$, in the resulting coefficient of the new $p$th order differential operator $\Alt\bigl( 1\cdot\dd_x^p,\ldots,x^{2p-1}\cdot\dd_x^p\bigr)$.
Having found the positive integer values $\consta(p\leqslant 20)$, we discovered that this integer sequence $\{\consta(p)\}_{p\in\BBN}$ is new in the \textit{On-Line Encyclopedia of Integer Sequences} (\cite{oeisA392714}).
We explore its properties in~\cite{Method,OddEven} and subsequent work (see~\cite{cpWron}). For instance, 
the prime factorisation of~$\consta(p=14)$ gave us a previously unknown, now confirmed \emph{prime} number of 206 decimal digits (\cite{Method}).

We expect that an explicit description of the sets $\Phi_p$ will unravel more properties of the sequence $\{\consta(p)\}_{p\in\BBN}$ (and will speed up the calculation of its terms beyond the current threshold $p=20$).
In particular, we hope that the technique which we employ in this paper will prove --\,for prime numbers~$p$\,-- the conjectured equality $v_p\bigl(\consta(p)\bigr)=p-1$ for the $p$-adic valuations of the positive integers~$\consta(p)$; this will strengthen the lower bound $v_p\bigl(\consta(p)\bigr) \geqslant p-1$.

\medskip\noindent\textbf{Acknowledgements.}\quad
The authors thank R.\,Buring, M.\,Kon\-tse\-vich, and O.\,Zaboronski for fruitful collaboration and discussions.
The first author thanks the Organisers of the 36th International colloquium on group theoretical methods in Physics (13--17 July 2026 in Valladolid, Spain) for a partial waiver of the registration fee; the travel to Valladolid was supported in part by the Bernoulli Institute (Groningen, NL) via project `Graduate programme'.
The second author thanks the audience of Mathematics seminar (held on 13 May 2026) and 13th symposium on Combinatorics and Arithmetic for Physics (CAP13, held on 26--27 November 2026) at the IH\'ES for stimulating remarks and insights; the travel to Bures\/-\/sur\/-\/Yvette was supported in part by the IH\'ES and by the Bernoulli Institute via project `Algebra'.

\paragraph{Author contributions}
Both authors contributed equally to this work.

\paragraph{Funding information}
No funding for the submitted work was obtained from any organisation.

\end{document}